\documentclass{amsart}
\usepackage{amsfonts,amsmath,amssymb,eucal}
\usepackage[all]{xy}
\usepackage{hyperref}
\usepackage[T1]{fontenc}
\usepackage[utf8]{inputenc}

\begin{document}
\newcommand{\M}{{\mathbb M}}
\newcommand{\bP}{{\mathbb P}}
\newcommand{\bH}{{\mathbb H}}
\newcommand{\Kerr}{{ \rm Kerr}}
\newcommand{\bx} {{\bf x}}
\newcommand{\bj}{{ \sf j}}
\newcommand{\vol}{{\rm vol}}
\newcommand{\w}{{\boldsymbol{\varepsilon}}}
\newcommand{\R}{{\mathbb R}}
\newcommand{\C}{{\mathbb C}}
\newcommand{\Z}{{\mathbb Z}}
\newcommand{\half}{{\textstyle{\frac{1}{2}}}}
\newcommand{\tk}{{\theta}}
\newcommand{\bX}{{\bf X}}
\newcommand{\D}{{\mathbb D}}
\newcommand{\T}{{\mathbb T}}
\newcommand{\Spin}{{\rm Spin}}
\newcommand{\SU}{{\rm SU}}
\newcommand{\U}{{\mathbb U}}
\newcommand{\E}{{\rm E}}
\newcommand{\kerr}{{\rm Kerr}}
\newcommand{\cayley}{{\sf C}}
\newcommand{\light}{{\lambda}}
\newcommand{\vw}{{\varpi}}
\newcommand{\sla}{{\star \lambda}}
\newcommand{\bl}{{\boldsymbol{\lambda}}}
\newcommand{\nun}{{\widetilde{\rm unknot}}}
\newcommand{\ob}{{\mathfrak o}}
\newcommand{\vk}{{\kappa}}
\newcommand{\arcsinh}{{\rm arcsinh}}
\newcommand{\pt}{{\rm pt}}
\newcommand{\ab}{{\rm ab}}
\newcommand{\lcs}{{\sf lcs}}
\newcommand{\balpha}{{\boldsymbol{\alpha}}}
\newcommand{\bbeta}{{\boldsymbol{\beta}}}
\newcommand{\bgamma}{{\boldsymbol{\gamma}}}
\newcommand{\ks}{{\sf ks}}
\newcommand{\Aa}{{\mathfrak a}}
\newcommand{\wU}{{\mathbb U^\ab}}
\newcommand{\vphi}{{\boldsymbol {\varpi}}}
\newcommand{\ve}{{\boldsymbol{\varepsilon}}}

\newcommand{\ie}{\textit{i}.\textit{e}.\,}
\newcommand{\eg}{\textit{e}.\textit{g}.\,}
\newcommand{\cf}{{\textit{cf} \,}}
\newcommand{\que}{{\mathord{?}}}

\parindent=0pt
\parskip=6pt

\title{A locally conformally symplectic structure on fast Kerr space-time}

\author{Jack Morava}

\address{Department of Mathematics, The Johns Hopkins University, Baltimore, Maryland 21218}

\email{jmorava1@jhu.edu}

% 83C15}

\begin{abstract}{The Kerr model for space-time around a rotating black hole is a classic solution of the vacuum equations of general relativity \cite{15,20}, a fundamental example in the study of angular momentum \cite{21}.
% a cosmological whirlpool, vortex, tornado made of no matter at all, not even a photon: the rushing % darkness of Pindar. [H\"olderlin] This note records an attempted solution to a textbook exercise in % its study, suggested by \cite{3}(Prop 3.1).}\end{abstract}\bigskip
In classical terms \cite{2} the Kerr solution defines a completely integrable system, but as thus presented it is fundamentally a birational object. Here we construct a version of Kerr space-time as a locally conformally symplectic double branched cover of the Penrose compactification of Minkowski space, and argue the interest of  the cobordism theory of 3D contact manifolds under 4D $\lcs$ cobordism.}\end{abstract}\bigskip

\maketitle \bigskip

{\bf Introduction} \bigskip

In classical general relativity, Kerr coordinates 
\[
\ks : \Kerr := \R \times (\R^3_+ -  \T_a)  \to \M 
\]
map a Ricci-flat Lorentz metric on the complement of a singular locus $\R \times (T_a = \{ (x,y,0) \:|\: x^2 + y^2 = a^2 \}) \subset \R^{1,3}$ to Minkowski space, generalizing the Schwarzschild metric on $\M - \R \times (0,0,0)$. Interpreting points of $\M$ as $2 \times 2$ self-adjoint matrices defines the Cayley transform 
 \[
 \cayley : \M \ni X \mapsto \frac{X - i}{X + i} \in \U(2)
 \]
 with values in the compact group of $2 \times 2$ unitary matrices \cite{14}. Acting on expert advice, we define locally conformally symplectic structures on these manifolds in terms of differential forms
 \[
 \omega_\Kerr \in \Omega^2 \Kerr, \; \omega_{\U(2)} \in \Omega^2 \U(2) 
 \]
(\S 2.2) and calculate $\omega_\Kerr$ in terms of $(\cayley \ks)^*\omega_{\U(2)}$. We then argue that when $a = 2$, generalizing Kerr coordinates by interpreting $\phi$ as taking values in $\R$ rather than $\R/\Z$ leads to a commutative covering diagram 
\[
\xymatrix{
\Kerr^+ \cong \R^2 \times \T \times \R \ar[d] \ar[r]^-{1 \otimes [2]} & \R^2 \times \T \times \R \ar[d] \\
\Kerr^\ab \ar[r]^-{\cayley \ks^\ab} & \wU(2) }
\]
presenting a lift of $\cayley \ks$ as a branched double cover. \newpage

{\bf \S 1 Some coordinate atlases} \bigskip

{\bf \S 1.1}  Euler angles for $\U(2)$

We first recall the Maurer-Cartan forms on $\U(2)$ and its universal cover:

Let $\bH = \R \langle \bf {i, \: j \: k} \rangle$ be the division algebra of quaternions, with
\[
q \mapsto (|q|,|q|^{-1/2}q) : \bH^\times \cong \R_{>0} \times \SU(2) 
\]
and define Euler angles or coordinates $\alpha,\beta,\gamma$ for $\SU(2)$ by 
\[
\xymatrix{
(\R/\Z)^3 \ar[d]^\cong \ar@{.>}[dr]^\varepsilon \\
\T^3 \ar[r] ^{\bf e} &  \SU(2) \;,}
\]
\[ 
{\bf e}(\alpha,\beta,\gamma) = e(\gamma{\bf k}) \cdot e(\beta{\bf j}) \cdot e(\alpha{\bf i})
\]
with $e(x) = \exp(2\pi x)$. It then follows formally that 
\[
\balpha = d\alpha + \cos \beta \:  d\gamma
\]
\[
\bbeta = - \sin \alpha \: d \beta + \cos \alpha \sin \beta  \: d \gamma 
\]
\[
\bgamma = \cos \alpha \: d \beta + \sin \alpha \sin \beta \: d \gamma  \in \Omega^1 \SU(2)
\]
satisfy 
\[
d\balpha = \bbeta \wedge \bgamma, \; d\bbeta = \bgamma \wedge \balpha, \; d\bgamma = \balpha \wedge \bbeta \in \Omega^2 \SU(2) \;.
\]
\cite{26,27}, see [Lie algebra cohomology] \dots \bigskip
% See for example {\tt http://compalg.jinr.ru/CAGroup/Palii/GeomSU23.pdf} 

We regard 
\[
\xymatrix{
{\SU(2) \subset \U(2)  = \T \times_{\pm 1} \SU(2) \ni  
\left[\begin{array}{cc} 
                        u & v \\
- \tau \bar{v} & \tau \bar{u}  \end{array}\right]} \ar[r]^-{\rm det} & \tau \in \T} 
\]
(with $|u|^2 + |v|^2 = 1$) as the mapping torus of an endomorphism of $\SU(2)$, a very basic example of a Milnor fibration. 
The fiber product
\[
\xymatrix{
\U^\ab(2) \cong \R \times \SU(2) \ar@{.>}[d] \ar@{.>}[r]^-t & \R \ar[d]^{\bf e} \\
\U(2) \ar@{.>}[r]^-\tau & \T \;,}
\]
(with $\tau = \exp(it) = {\bf e}(t)$) defines a lift 
\[
\Omega^1\U(2) \ni i d\tau \mapsto t^{-1} dt \in \Omega^1\U^\ab(2) 
\]
and thus a set $\balpha, \bbeta, \bgamma, dt$ of translation-invariant elements in $\Omega^1 \U^\ab(2)$.
% cf ${\rm Spin}^c(3)$. 
\bigskip

The Cayley transform defines a homomorphism
\[
\cayley^* : \Omega^*\U^\ab(2) \to \Omega^*\M
\]
of DG algebras of differential forms, perhaps useful as characteristic classes for geometric questions in Lorentz/Minkowski topology. \medskip

{\bf Exercise} 
\[
\half \: {\rm Trace} \; \cayley(t^{-1}\bX) \sim 1 +  \sum_{n \geq 1} {\rm Trace} \; \bX^{-n} \cdot (-it)^n
\]
is something like a cumulant \dots \bigskip

{\bf 1.2 $\lcs$ four-manifolds}

In 1969 Brandon Carter \cite{2} showed that the system of geodesics in Kerr's Lorentz Ricci flat Einstein space-time vacuum is completely integrable: which he then (\S A3) integrated. Recent language interprets this as part of a Poisson structure on a locally conformally symplectic (or $\lcs$) manifold \cite{6, 10, 16, 19}. Such structures have a rich theory of deformation quantization. 

I am deeply indebted to conversations in Summer 2022 with YM Eliashberg, when he explained me that an $\lcs$ structure on a four-manifold $X$ is (roughly) a two-form $\omega$ such that $\omega \wedge \omega $ is (as in symplectic geometry) a volume form, but is more general in that rather than its being closed, we require
\[
d\omega = \eta \wedge \omega
\]
with the (closed) Lee one-form \cite{11} $\eta$ classified by a flat line bundle
\[
[\eta] \in {\rm Pic}_\R(X) \cong {\rm Hom}(\pi_1X,\R^\times) \cong H^1(X,\R^\times) \;. 
\]
Such structures have interesting Kolmogorov-Arnol'd-Moser stability properties \cite{1}, and lead thus to toric manifolds\dots
\bigskip

{\bf Example} \cite{5}(\S 3.2),\cite{3}(Prop 3.1) : {\it If $X = \U(2)$ then 
\[
\omega_{\U(2)} = d \balpha + \balpha \wedge \eta_{\U(2)} \ni \Omega^2 \U(2)
\]
with $\eta_{\U(2)} = id\tau$ and 
\[
\omega_{\U(2)} \wedge \omega_{\U(2)} = 2 \balpha \wedge \bbeta \wedge \bgamma \wedge i d\tau \in \Omega^4 \U(2)
\]
defines an $\lcs$ structure on $\U(2)$, which lifts to the universal cover $\wU(2)$.} \bigskip

{\bf 1.3} Kerr-Schild Cartesian coordinates \bigskip

{\bf 1.3.1}  Some conventions : It is useful to distinguish the Kerr {\it coordinate system} from the Kerr metric itself; this note is more about the former (an atlas for the manifold underlying the model, and variants) than the latter; see the appendix below for some details.

Composing Kerr-Schild Cartesian coordinates \cite{20}(eq 34, hereafter abbreviated as V eq ...) with the Cayley transform defines 
\[
\cayley \ks : \R^2 \times S^2  \ni (u_\pm,\theta,\phi) \mapsto (t,\bx) = (t,x,y,z) \mapsto (x_0,\dots,x_3)
\] 
\[
= x_* \to \left[\begin{array}{cc}
     x_0 + x_3 & x_1 + ix_2 \\
     x_1 - ix_2 & x_0 - x_3 \end{array}\right]
     =  \bX \in {\rm M}_2(\C) \to \cayley(\bX) \in \U(2) \;.
\]
Here $u_\pm = t \pm r$,  with $\; t,r \in \R$; note that $r$ is distinct from $|\bx| = (x^2 + y^2 + z^2)^{1/2}$. We work in Planck units with $m=1$, and in practice mostly with $t$ and $r$. In Minkowski space the advanced and retarded \cite{22} parameters $t \pm |\bx|$ are the eigenvalues of $\bX$.

We have
\[
x + iy = (r - ia) \: \sin \theta  \:e^{i\phi}, \; x_3 = z = r \cos \theta \;,
\]
so [V eq 33]
\[
x^2 +  y^2 + z^2 = a^2(1 - r^{-2}z^2) \;,
\]
which makes it useful to recall the \bigskip

{\bf lemma} \cite{18}(\S 15.24),\cite{23, 26}
\[
 J_\pm(x) = \half (x \pm x^{-1})
\]
{\it defines an invertible map $J_+ : \R^\times \to \R^\times$ and a double cover $J_- : \R^\times \to \R_{>0}$ satisfying $iJ_-(x) = J_+(ix)$. Moreover,
\[
J_-(x + \surd (x^2 + 1)) = x 
\]
and hence 
\[
J_-^{-1}(x) := x + \surd (x^2 + 1) > 0 \;.
\] } \medskip

{\bf 1.3.2} More conventions: Ignoring singularities for the moment, let us define
\[
r^2/az = s \: \exp \varkappa , \; \; \; \ob(\bx) := \frac{|\bx|^2 - a^2}{2az} , 
\]
where $\pm 1 \ni s =$ sign of $z \in \R$.  Both $\ob$ \cite{24} and $\varkappa$ are invariant under rescalings of length, and are `pure numbers' in the sense of physics. See [Further digressions] below. 
\bigskip

{\bf Proposition} 
\[
\varkappa = \arcsinh \: s \ob, \; \; r =  a \cos \theta \cdot sJ_-^{-1}(s \ob) \;.
\]

{\it When $\ob \neq 0$ the level surfaces $|\bx|^2 = (a + \ob z)^2 - (\ob z)^2$ are conics ($z \sim r^2$) 
but if $r = 0$ then $x + iy = - ia \sin \theta \: e^{i\phi}$ parametrizes the disk $x^2 + y^2 \leq a^2, \; z = 0$. On the other hand, if $|\bx|> 0$ then $r \to \pm \infty$ as $z \to \pm 0$.}\bigskip

{\bf 1.3.3} Following [Hawking $\&$ Ellis Fig 27], Kerr-Schild coordinates
\[
\ks : \Kerr := \R \times (\R^3_+ - \T_a)  \to \M \;,
\]
(where $\T_a$ is a circle of radius $a$, or, alternatively, a disk rotating with angular acceleration $a>1$) map Kerr space-time to a constant unknot in Minkowski space. 
 
A beautiful theorem of Stallings and Neuwirth asserts that if the commutator subgroup of the fundamental group $\pi_1(S^3 - k)$ of a knot complement is finitely generated, then $S^3 - k$ is a bundle over $S^1$ with fiber a punctured surface  $\Sigma(k) - \pt$. Its universal abelian (or Alexander \cite{17}(Ch 7)) cover $(S^3 - k)^\ab$ then splits \cite{12} as a product $\R \times (\Sigma(k) - \pt)$, with deck-transformation
\[
\vphi : (\tau,x) \mapsto (\tau + 1,\varpi(x))
\]
defined by a basepoint-preserving diffeomorphism of $\Sigma$. \newpage

{\bf Definition}
\[
\R \times_\Z (S^3 - k)^\ab \;=\;  \mathfrak{M}(k) \;,
\]
with the standard Lorentz structure, is the {\it metaverse}$\textsuperscript{\textcopyright}$ indexed by $k$.\bigskip

[Stereographic projection embeds $\R^{1,3}$ in $\R \times S^3 \cong \R \times \R^3_+$, with $\{+ = \infty\}$ as basepoint. Note that for a knot $k$, $\R^3_+ - k$ (unlike $\R^3 - k$) is acyclic: it has no higher homotopy groups, which simplifies some covering space issues: a cyclic cover of $S^3$ branched over the circle is diffeomorphic to $S^3$, \cite{17}(Ch V \S C, Ch X ex 5).]

In this quotient, $1 \in \Z$ acts as the unit translation of $\R$, and as $\vphi$ on the second factor. This provides a large supply of potentially interesting Lorentz four-manifolds; the unknot, with
\[
\pi_1(S^3 - S^1) \cong \Z
\]
and $\Sigma \cong S^2$, is an example. The Kerr-Schild map
\[
\ks : \R^2 \times (S^2 \supset (\R^2  = S^2 - \infty)) \to \R \times (\R^3_+ - \T_a)
\]
(as defined in \S 1.3.1) extends/lifts naturally to the infinite cyclic cover of the unknot: in particular, the Kerr metric is well-defined for $\phi \in \R$, not just for $\phi \in \T \cong \R/\Z$. We now have a fiber product
\[
\xymatrix{
{\R^2 \times \R^2 := \Kerr^+}   \ar@{.>}[d]  \ar@{.>}[r]  & \R \times \R^3 \ar[d]^\ve \\
{\Kerr^\ab \cong \R^2 \times S^2}  \ar[r]  & \U^\ab(2) }
\]
defining a DGA morphism $\cayley \ks^* : \Omega^* \wU(2) \to \Omega^* \Kerr^+$. \bigskip \bigskip

{\bf \S 2 Hodge theory on Kerr space-time} \cite{20}(\S 3-4) \bigskip

{\bf 2.1} Following M Visser's insightful presentation closely, we collect some basic material, postponing some relevant details:\bigskip

The Kerr solution is characterized by a field $\bl = (1,0,0, a \sin^2 \theta)$ of null vectors; the one-form [V eq 14]
\[
\lambda = \sum \lambda_i dx_i = du + a \sin^2 \theta \: d\phi \in \Omega^1 \Kerr
\]
is such that 
 \[
d\lambda = a \sin 2 \theta \; d\theta \wedge d\phi \in \Omega^2 \Kerr \;,
\]
\[
\lambda \wedge d\lambda = a \sin 2 \theta \: du \wedge d\theta \wedge d\phi = 2 \nu \cdot \sla \in \Omega^3 \Kerr \;,
\]
and
\[
\star 1 = d\vol_\Kerr =  - (r^2 + a^2 \cos^2\theta)^2 \sin^2 \theta \cdot dr \wedge d \theta \wedge d\phi \wedge dt \in \Omega^4 \Kerr \;,
\]
where $\star$ is the Lorentz Hodge operator [V eq 29, wiki]. \newpage

If we write
\[
\nu = - \frac{a \cos \theta}{r^2 + a^2 \cos^2 \theta} = 2r \: {\rm sech} \; \varkappa  = (1 + \ob^2)^{-1/2} 
\]

[V eq 31] for the twist  denoted $\omega$ in \cite{20}, then 
\[
d\vol_\Kerr = - \nu^{-2} \cdot a^2 \cos^2 \theta \sin^2 \theta \cdot dr \wedge d \theta \wedge d\phi \wedge dt 
\]
and 
\[
\lambda \wedge d \lambda \wedge dt = - 4a^{-1} \csc 2 \theta \cdot \nu^2 \: d\vol_\Kerr \;.
\]
\bigskip

{\bf Corollary} $\lambda \wedge d \lambda$ is nonvanishing on $t =$ constant hypersurfaces, whereupon $\lambda$ (reminiscent of $dx + ydz$) {\bf defines a contact form}.
\bigskip

{\bf 2.3} Now 
\[
d(t\lambda) = - \lambda dt  + t d\lambda = - (du + a \sin^2\theta) \wedge dt + t \cdot a \sin 2 \theta d\theta \wedge d\phi = -dr \wedge dt + d \Theta \wedge d\phi 
\]
with $\Theta = at\sin^2 \theta$, so 
\[
d(t\lambda) \wedge d(t\lambda) = - 2 dr \wedge dt \wedge d\Theta \wedge d\phi = - 2a \sin 2 \theta \: dr \wedge d \theta \wedge d\phi \wedge dt 
\]
\[
= \frac{-2a \sin 2 \theta}{- \nu^{-2} a^2 \sin^2 \theta cos^2 \theta} \: d \vol_\Kerr 
\]
and hence 
\[
t^{-1}d(t\lambda) \wedge t^{-1} d(t\lambda) = 8a^{-1} \sin 2\theta \cdot (t^{-1} \nu)^2 \; d \vol_\Kerr \;. 
\]
\bigskip

{\bf Definition} If
\[
\omega_{\Kerr^+}  := t^{-1}d(t\lambda) = d\lambda - \lambda \wedge i t^{-1}dt = d\lambda + \lambda \wedge \eta  \;.
\]
then 
\[
\omega_{\Kerr^+} \wedge \omega_{\Kerr^+}  =t^{-1}d(t\lambda) \wedge t^{-1} d(t\lambda) = 8a^{-1} \sin 2\theta \cdot (t^{-1} \nu)^2 \; d \vol_\Kerr \;. 
\]
while 
\[
d\omega_{\Kerr^+} = d(d\lambda + \lambda \wedge i d\tau) = \omega_{\Kerr^+} \wedge \eta \;.
\]
\bigskip

We thus have $\lcs$ structures 
\[
\omega_{\U^\ab(2)} = d\alpha + \alpha \wedge \eta
\]
and 
\[
\omega_{\Kerr^+} = d\lambda + \lambda \wedge \eta
\]
with $\eta =  i d\tau = t^{-1}dt$ on both $\wU(2)$ and $\Kerr^+$, and the graph
\[
\cayley \ks : \{ (t,u,\theta,\phi) = (t,\alpha,\beta,\gamma)\} \subset \Kerr^+ \times (\R \times \R^3) :
 \Kerr^+  \to \R \times \R^3
\]
 which defines a map between them. \bigskip
% ? maybe \Spin^c(3)

{\bf Proposition} {\it Writing
\[
[a]  \; :=  \;  \left| \begin{array}{ccc}
                                        1  &  0  &  a/2 \\
                                        0  &  2  &   0   \\
                                        0  &  0  &  -a/2  \end{array} \right| \in {\rm Gl}_3(\R) \;,
\] 
we have 
\[
\lambda = [a] (\cayley \ks^*(\balpha)) \in \Omega^1 \Kerr^+ 
\]
and 
\[
\omega_{\Kerr^+} = [a] (\cayley \ks^* \omega_{\wU(2)})  \in \Omega^2 \Kerr^+ \;.
\] }
\bigskip

{\bf Proof:} 
\[
\lambda = du + a \sin^2 \theta d\phi  = d(u + a/2 \; \phi) + \cos 2 \theta  \cdot d(- a/2 \: \phi) \; ;
\]
that is, 
\[
\lambda =  [a] \cayley \ks^*(\alpha + \cos \beta \: d \gamma) \;  \;.
\]
Similarly, by linearity
\[
\omega_{\Kerr^+}  = (d + \eta) \lambda  = [a] (d + \eta) \cayley \ks^* (\balpha) = [a]  \cayley \ks^* \omega_{\wU(2)}  \; \;.
\]
\bigskip

Google says $[a]$ has minimal polynomial $2T^3 + (a - 6)T^2  + (4 - 3a)T + a$, and we have
\[
[2]  \;\; = \; \; \left| \begin{array}{ccc}
                                        1  &  0  &  1 \\
                                        0  &  2  &   0   \\
                                        0  &  0  &  - 1 \end{array} \right| \;.
\]
\bigskip

{\bf Corollary} The commutative diagram 
\[
\xymatrix{
{ } & \R^2 \times \T \times \R \ar[dl] \ar[d] \ar[r]^-{1 \otimes [2]} & \R^2 \times \T \times \R \ar[d] \ar[dr] \\
\R^2 \times S^2 \ar[r] & \Kerr^\ab \ar[r]^-{\cayley \ks} & \wU(2) & \R \times \T^3 \ar[l] }
\]
defines a completion of Kerr space-time as a double cover of $\wU(2)$. \bigskip

This suggests comparison of the associated Poisson structures \cite{6} under maps such as $1 \times [2]$. See also the horizons of \cite{28}. \bigskip

{\bf \S 3} \; (in progress) The cobordism category
\[
{{\sf IV}^{\rm lcs}}_{\rm ctc}
\]
has closed compact oriented {\bf contact} three-manifolds as objects, and compatible locally conformally symplectic four-manifolds as cobordisms between them. More precisely, a morphism
\[
X \supset \partial X = Y^{\rm op}_- \cup Y_+ : Y_ - \to  Y_+
\]
is a morphism of such structures modulo conformal equivalence as in \cite{5}|(\S 3, 6.3), \ie
\[
{\rm Mor}_{\sf IV}(Y,Y') :=   \coprod_{\pi_0 \ni X} B{\sf HR}(X,\partial X)
\]
with
\[
{\sf HR} := {\rm Aut}_\lcs /{\rm conf} ;.
\]
This is a topological category \cite{9} with an interesting Abel - Jacobi completion in terms of Calabi fluxes and Picard categories of line bundles classified by two-dimensional Lichnerowicz cohomology \cite{5}(\S 6) \dots \bigskip

{\bf \S 4 Remarks}

Relations between coverings and compactifications can be complicated. The diagram
\[
\xymatrix{
(S^3 - S^1)^\ab \cong \R \times (S^2 - \pt) \ar[r] &  \R \times_{2\Z} (S^2 - \pt) \ar[d] \\
{ } & S^3 - S^1 }
\]
(\cf O'Neill \cite{15}[Fig 3.20]) presents a partial compactification of the universal cover of Kerr space-time as a homotopy suspension
\[ 
\R^2 \times \bigvee_\Z S^2 \simeq S^1 \wedge (S^2 \vee S^1)^\ab
\]
of Milnor's string \cite{13} of balloons.
% , bubbles, pearls of Indra, the celestial two-sphere, the night sky. 
A model is perhaps the universal cover of a torus $\T^2/\T \times 0$ with a pinched cycle, a string of conformal blowdowns of Minkowski space-times at infinity resembling the bubble-chains seen in particle detectors. See further below re the hypervolume $\textstyle{\frac{5}{4}}\pi^3$ of a unit Friedman aeon. \bigskip

{\bf For reference}: the Kerr metric itself [V eq 43, 44]
\[
g_{ij} = \eta_{ij} + \frac{2r^3}{r^4 + a^2 z^2} \; \lambda_i \lambda_j \;, \; \eta_{00} = -1 \;
\]
is a rank one deformation \cite{12} of the Minkowski metric, defined by a field 
\[
\bl(r,\bx)  = (1, \frac{rx + ay}{r^2 + a^2}, \frac{ ry - a x}{r^2 + a^2}, r^{-1}z)
\]
($\bl = (1,0,0, a \sin^2 \phi)$ in the original polar coordinates), of null geodesics with gradient [V eq 22]
\[
\nabla \cdot \bl = \frac{2r}{r^2 + a^2 \cos^2 \theta} = - 2s(1 + \ob^2)^{-1/2} \exp \varkappa
\]
\cf Geroch \dots. This metric is asymptotically Minkowskian, which suggests compactifying its underlying manifold as a pushout of a compactification of $\M$, but it may be useful to pass to a suitably universal cover first. \newpage

{\bf Further digressions} 
 
 re \S 1.3.2 : Physical quantities are often measured in units, and are thus subject to rescaling; the parameter $a$ in the Kerr model, for example, scales like a length. A nontrivial scaling law suggests interpreting such a quantity as a section of a line bundle; see \cite{29}(\S 4.2) for the Planck length $(Gh)^{-1/2} \sim 10^{-17}$ Higgs.
 
 re \S 4 : The Cayley transform defines a compactification of Minkowski space in terms of a conformal factor which vanishes on the two-sphere at infinity in $\U(2)$ defined by unitaries with eigenvalues $\pm 1$. Collapsing the coverings of this bubble in $\wU(2)$ to points defines a string of `aeons', familiar in conformal variants of classical GR. \bigskip
 
 % at 6:30 PM 13.12.022 I learned that there is a well-established and very interesting notion of an acoustic metric \cite{28}.

\bibliographystyle{amsplain}

\begin{thebibliography}{99}

\bibitem[1]{1} G Bande, D Kotschick, Moser stability for locally conformally symplectic structures, Proc. Amer. Math. Soc. 137 (2009) 2419 -- 2424
 
\bibitem[2]{2} B Carter, Global structure of the Kerr family of gravitational fields, \url{https://journals.aps.org/pr/pdf/10.1103/PhysRev.174.1559} \S 3A

\bibitem[3]{3} YM Eliashberg, On symplectic manifolds with some contact properties, J. Differential Geom. 33 (1991) p 233 –- 238

\bibitem[4]{4} G Gibbons, S Hawking, Kinks and topology change, Phys. Rev. Lett. 69 (1992) 1719 -– 1721

\bibitem[5]{5} S Haller, T Rybicki,  On the group of diffeomorphisms preserving a locally conformal symplectic structure. Ann. Global Anal. Geom. 17 (1999) 475–502

\bibitem[6]{6} -----, -----, Integrability of the Poisson algebra on a locally conformal symplectic manifold, 19th Winter School, Rend. Circ. Mat. Palermo (2) Suppl. No. 63 (2000) 89 -–  96

\bibitem[7]{7} C Heinicke, F Hehl, Schwarzschild and Kerr solutions \dots, \url{https://arxiv.org/abs/1503.02172}

\bibitem[8]{8} W Israel, Source of the Kerr metric, Phys Rev D2 641 (1970)

\bibitem[9]{9} J Giansiracusa, The stable mapping class group of simply connected 4-manifolds. J. Reine Angew. Math. 617 (2008) 215 –- 235, \url{https://arxiv.org/abs/math/0510599}

\bibitem[10]{10} M Kontsevich, Deformation quantization of Poisson manifolds, Lett. Math. Phys. 66 (2003) 157 -– 216, \url{https://arxiv.org/abs/q-alg/9709040}

\bibitem[11]{11} Hwa-Chung Lee, A kind of even-dimensional differential geometry and its application to exterior calculus. Amer. J. Math. 65 (1943), 433 – 438

\bibitem[12]{12} M Mars, Spacetime Ehlers group: transformation law for the Weyl tensor, \url{https://arxiv.org/abs/gr-qc/0101020}

\bibitem[13]{13} J Milnor, Infinite cyclic coverings, in 1968 Conference on the Topology of Manifolds (Michigan State Univ. 1967) 115 – 133 Prindle, Weber \& Schmidt, Boston

\bibitem[14]{14} J Morava, At the boundary of Minkowski space,  \url{https://arxiv.org/abs/2111.08053}

\bibitem[15]{15} B O'Neill,{\it The geometry of Kerr black holes}, A K Peters, Ltd., Wellesley, MA, 1995.

\bibitem[16]{16} A Otiman,  Locally conformally symplectic bundles, J. Symplectic Geom. 16 (2018) 1377 -– 1408, \url{arXiv:1510.02770}

\bibitem[17]{17} D Rolfsen, {\it Knots and links}, Mathematics Lecture Series 7, Publish or Perish, Inc., Berkeley 1976

\bibitem[18]{18} RA Silverman, {\it Complex analysis with applications}, Dover 1984

\bibitem[19]{19} I Vaisman, On the geometric quantization of Poisson manifolds. J. Math. Phys. 32 (1991) 3339–3345

\bibitem[20]{20} M Visser, The Kerr spacetime: A brief introduction \url{https://arxiv.org/abs/0706.0622}

\bibitem[21]{21} Po-Ning Chen, Mu-Tao Wang, Ye-Kai Wang,  Shing-Tung Yau, Supertranslation invariance of angular momentum, Adv. Theor. Math. Phys. 25 (2021) 777 -– 789

\bigskip

\bibitem[22]{22} \url{https://en.wikipedia.org/wiki/Kruskal%E2%80%93Szekeres_coordinates}

\bibitem[23]{23} \url{https://en.wikipedia.org/wiki/Inverse_hyperbolic_functions}

\bibitem[24]{24} \url{https://en.wikipedia.org/wiki/Oblate_spheroidal_coordinates}

% \bibitem[25]{25} \url{https://en.wikipedia.org/wiki/Electroweak_interaction}

\bibitem[25]{25} \url{https://en.wikipedia.org/wiki/Joukowsky_transform}

\bibitem[26]{26} \url{https://en.wikipedia.org/wiki/Euler%27s_equations\_(rigid\_body\_dynamics)}

\bibitem[27]{27} \url{https://math.stackexchange.com/questions/3706481\_explicit\_computation\_of\_a\_3-form\_on\_su2}

% \bibitem[27]{27}\url{https://en.wikipedia.org/wiki/Milne_model} 

\bibitem[28]{28} \url{https://en.wikipedia.org/wiki/Acoustic_metric}

\bibitem[29]{29} JM, Conformal invariants of Minkowski space, Proc. Amer. Math. Soc. 95 (1985) 565 – 570 [ MR0810164 ], \;  On gauge theories of mass, J. Geom. Phys. 62 (2012) 1262 – 1272 [MR2901859], \url{https://arxiv.org/abs/1001.0965},

\end{thebibliography}

\bigskip

\end{document}